\documentclass[10pt]{article}
\usepackage{amsmath}
\usepackage{graphicx}
\usepackage{hyperref}
\usepackage[T1]{fontenc}
\usepackage{helvet}

\begin{document}           
\title{The Amazing Journey of Lars Ahlfors' Fields Medal}
\author{Frank Wang, LaGuardia Community College,\\
  The City University of New York\\
  \url{fwang@lagcc.cuny.edu}}
\date{}

\maketitle

Lars Ahlfors was born in Finland in 1907.  In 1936, he received, along
with Jesse Douglas, one of the first two Fields Medals.  At the time,
he was a professor at the University of Helsinki, and a visiting
professor, from 1935 to 1938, at Harvard University.  In the 1997 film
\textit{Good Will Hunting}, the psychology professor portrayed by
Robin Williams explained to his community college students that the
Fields Medal is ``like the Nobel Prize for math, except they only give
it out once every four years.''  Furthermore, the recipient must be
under forty years of age.  The Fields Medal is widely considered to be
the most prestigious mathematics award.

Elevated by the recognition in the mathematical world, the Ahlfors
couple returned to Finland from Harvard, and had their first child
Cynthia in 1938.  Their happiness did not last long.  Nazi Germany
invaded Poland on September 1, 1939, and Europe descended into chaos.
The second daughter Vanessa was born a month after the Second World
War broke out.  In 1944, Ahlfors was in a desperate situation to leave
Finland.  His account of the escape was quoted in his obituary written
by Steven Krantz published in \textit{Notices of the American
  Mathematical Society}, based on Donald Albers' 1991 interview of
Ahlfors published posthumously in the \textit{College Mathematics
  Journal} (\textit{CMJ}).
\begin{quote}
I can give one very definite benefit [of winning a Fields Medal].
When I was able to leave Finland to go to Sweden, I was not allowed to
take more than 10 crowns with me, and I wanted to take a train to
where my wife was waiting for me.  So what did I do? I smuggled out my
Fields Medal, and I pawned it in the pawn shop and got enough money.
I had no other way, no other way at all.  And I'm sure that it's the
only Fields Medal that has been in a pawn shop.
\end{quote}

Krantz's story of the Fields Medal in a pawn shop was included in his
amusing book \textit{Mathematical Apocrypha, Stories and Anecdotes of
  Mathematicians and the Mathematical}.  I remember first reading it
after I bought the book during the Mathematical Association of America
MathFest held in Hartford, Connecticut in the summer of 2013.  I
chuckled occasionally on stories in the book on my Amtrak ride back to
New York, which caused my fellow passengers to raise their eyebrows.
But I also remember that a year earlier, on April 25, 2012, Ahllfors'
second daughter, now Vanessa Gruen and a longtime resident of Long
Island City, New York, where LaGuardia Community College is located,
came to our campus to give a math colloquium reminiscing about her
father.  She mentioned that she and her sister had donated the Fields
Medal to the University of Helsinki.  As I was unable to find the full
story of the first Fields Medal using internet search engines or
artificial intelligence tools at the time of this writing, I reached
out to Vanessa, and she kindly granted an interview.

Vanessa described that when she was a newborn, less than two months
old and was still in the hospital, Russia commenced the Winter War
against Finland.  Women and children were sent out of the country.
Her mother brought Cynthia and her to Sweden; she was so little and on
the boat she was mistaken by other passengers as Cynthia's doll.  The
Ahlfors family was forced to separate because Lars, as a young man,
was expected to stay to defend Finland.  When the fighting subdued,
the mother and daughters returned to Helsinki.  In the preface of
Ahlfors' \textit{Collected Papers}, which is a short autobiography, he
wrote the following.
\begin{quote}
Soon after the end of the Winter War, my family was able to return
home and resume a seemingly normal life.  Politics in Finland took an
unfortunate turn, however, and when Hitler attacked the Soviet Union
in 1941 Finland was his ally.  When the Russians were finally able to
repulse the attack, they could also intensify the war in Finland with
foreseeable results.  The Finish-Russian war ended with a separate
armistice in September 1944, whereupon Finland was able to expulse the
German troops stationed there.  The harsh terms of the armistice left
Finland in a very difficult position.
\end{quote}

In the \textit{CMJ} interview, Ahlfors said ``Being a Finnish citizen,
I was never in the military service, but I was still considered
eligible for special service.''  Vanessa explained that her father was
not taken as a soldier because he had a heart defect.  Instead, he was
sent to a room in a bomb shelter with some other mathematicians, and
they were supposed to code military messages.  But the phone never
rang and so the four mathematicians sat and played bridge the whole
time.

Vanessa said there were two Finish-Russian wars and asked me to look
it up; I did that.  The first was the above-mentioned Winter War, from
November 30, 1939 to March 13, 1940.  The second one, also called the
Continuation War, began on June 25, 1941 and ended on September 19,
1944.  In 1942, Ahlfors sent his two daughters to Sweden, to stay with
his wife Erna's sister in Kungsbacka.  In 1944, Erna and their young
son Christoffer (born in 1943) also left for Kungsbacka.  With an
offer from the University of Z{\"u}rich and worsening situation in
Finland, Ahlfors finally decided to emigrate permanently.  After
receiving clearance from the government, he left his homeland, with
his Fields Medal, in August 1944.

As told in his \textit{CMJ} interview, Ahlfors was permitted to take
only ten crowns out of the country, which did not cover the train fare
from Stockholm on the east coast of Sweden, to Kungsbacka, near
Gothenburg on the west coast of the country.  He left his Fields Medal
at a pawn shop for cash.  According to Vanessa, the mathematician Arne
Beurling offered her father a position at the Uppsala Institute near
Stockholm, and ``when he started getting his salary from Uppsala he
was able to take [the pawned Fields Medal] back.''  Ahlfors wrote ``I
am forever indebted to Arne Beurling, who showed what true friendship
can be.''

The Ahlfors family now faced the problem of getting from Sweden to
Switzerland.  When Finland signed an armistice with the USSR to end
the Continuation War, one of the conditions was to expel German forces
from its territory.  So in 1944, Finland was at war with Germany and
simultaneously, at least on paper, with the Allied nations including
the United Kingdom.  The British were willing to let the family pass
through the UK if an opportunity arose.  Ahlfors recalled that ``the
Swedes had organized some semiregular `stratospheric flights' on
moonless nights from Stockholm to Prestwick, Scotland.''  From
Vanessa's description, they had to wait for the phone to ring so that
they could know when to get on the plane; it was kept a secret.  While
the family was waiting in Sweden, a tragedy happened.  The young
Christoffer died in an accident.  It was understandably a very painful
and traumatic experience for the Ahlfors family, particularly for the
mother.

Ahlfors' account of their departure from Sweden is the following.
\begin{quote}
One day in March 1945 we were told to be ready to leave [Sweden],
weather permitting.  It is difficult to forget that flight.  The plane
was a reconditioned Flying Fortress, with perhaps a dozen passengers.
It was not pressurized, and breathing was accomplished by individual
oxygen masks.  Swim vests were worn by all.  Our children, ages 5 and
6, were quite capable of understanding the implications of danger.
\end{quote}

According to Vanessa, she was not sure if her oxygen mask fit properly
and didn't think it was working, but she was fine when they landed in
Scotland.  A train ride took the family from Glasgow to London.

Ahlfors recalled the continuing trip in his \textit{CMJ} interview.
\begin{quote}
After getting to England finally, we then had to cross the English
Channel to France.  But waiting for the boat-train took a long time.
We would take everything to the train station, where we sat on the
train, waiting for the boat.  At last they would come and tell us that
there was no boat that night, and we would go back to the hotel.  This
was repeated perhaps ten times, until finally one day the boat came
and the train truly left.  This was in early 1945.  London, of course,
was being bombed by the V-2s and the children were very scared.
\end{quote}

Vanessa said that the V-2 bombs were coming in during the day at any
time, and they could hear the big explosions.  But the father
explained to the two girls that it was only the tires on the
double-decker buses that were exploding, so as to pacify them.  Then
they proceeded to Paris, and from there took a train to Z{\"u}rich.  She
remembered that they arrived in a hotel in Z{\"u}rich, right before
Easter.  Her father had hidden a Swiss watch in a chandelier in the
hotel, which was the Easter gift for Erna.

Although Ahlfors said ``the Swiss legation was unforgettably helpful
in securing lodgings in a luxury hotel'' when they were in Paris,
which Vanessa identified as George Cinq, Ahlfors also said he ``cannot
honestly say that [he] was happy in Z{\"u}rich.  The postwar era was
not a good time for a stranger to take root in Switzerland.''  Vanessa
didn't think the country was very welcoming to foreigners.  From her
perspective, people were very stingy; they had all saved up coal.
Because the Ahlfors family had just arrived, they didn't know the
local custom and the apartment they stayed in was not adequately
heated.  The youngest daughter of the Ahlforses, Caroline, was born in
Z{\"u}rich in November 1945.  They put her in an orphanage, and she
remained there for three or four months until it became a little
warmer.  Vanessa remembered that when Caroline came back from the
orphanage, all of a sudden there was a big baby instead of the little
one that they had seen.

In 1946, Ahlfors accepted a professorship at Harvard University, and
stayed there until his retirement in 1977.  When Vanessa was in fifth
grade (around 1950), the family moved to Winchester, Massachusetts, a
suburb of Boston.  It was a spacious house, and Ahlfors had a home
office.  Based on Vanessa's memory, he would disappear and go to his
office every night after dinner.  Ahlfors was at the time writing the
textbook \textit{Complex Analysis}, which became tremendously popular
after its publication.  Ahlfors kept his Fields Medal, along with many
other medals, in his desk drawer.  Occasionally he would open the
drawer and show the girls all the medals.  Some of them were really
big, and the Fields Medal was a smaller one.  Ahlfors told his
daughters that ``this is the one that's worth a lot of money because
it's pure gold and others were not.''  Vanessa remembered going
through the medals in the drawer and she felt elated.

In the \textit{CMJ} interview, the daughter said ``Growing up, we
always were in awe of what father did.  We had no idea what it was.
We still don't know.  He tells us he's a philosopher.''  During our
interview Vanessa reiterated her feeling of her father being more a
philosopher than a mathematician; he never explained his math.  When
Vanessa was in high school and had math homework, she just went into
his office.  Ahlfors would say ``this is how you do it,
da-da-da-da-da-da, and he would do it.  Then I had my homework done.
He didn't try to teach it to me.  He just did the whole, I mean, it
was just easier to get it done.''

After finishing \textit{Complex Analysis}, Ahlfors continued to work
all the time, except Sundays.  He made everyone go to church.  His
family was Lutheran, like everybody else in Finland, but he had
converted to Catholicism and very much believed in that.  He was very
taken about the readings of Thomas Aquinas when he was around 18 or 19
years old, and he adhered to it.  The girls were reluctant.  Finally
when Vanessa was in college (New York University) and she was visiting
her parents in Massachusetts, an Irish priest in Boston said every
Catholic is obligated to vote for John Kennedy.  Ahlfors took great
offense at that, and they walked out of the church and never returned.
At the end, when he died, he was not at all Catholic.

Ahfors was a great pianist.  Vanessa said her father and all the
mathematicians that they knew at Harvard had a grand piano, and they
all played it.  Ahlfors would play in the afternoon or in the evening,
especially on Sundays.  He loved Bach and played all of his music.  He
also enjoyed playing Beethoven and Mozart.  And towards the end of his
life, Schubert became his favorite composer.

Ahlfors passed away in 1996, and his wife Erna in 2001.  Vanessa Gruen
and her sister Caroline Mouris inherited the Fields Medal.  (Cynthia
had died in 2000.)  They were passing it around from one coffee table
to the other, between Manhattan and Chatham in upstate New York.
Eventually they decided to donate the first Fields Medal won by Lars
Ahlfors to the University of Helsinki.  This is the quote from
Vanessa.
\begin{quote}
We were still in Chelsea in New York [a neighborhood in Manhattan].
We were on 7th Avenue, and my sister was living up in the Chatham area
in New York.  And it was just that we would switch [the Fields Medal]
back and forth, and we decided that was just ridiculous, that we
should just do something with it and give it something where it would
have some meaning.  And we thought that since he had done all his
mathematics and all this early work, and leading up to the Fields
Medal, it was all done in Helsinki.  It wasn't done at Harvard, and we
didn't want to give it to Harvard because they have so much stuff
anyway.  That would be more meaningful to give it to the University in
Helsinki where he had learned all this mathematics and where he was
when he was in Finland when he won the Fields Medal.  So we thought
that was appropriate.
\end{quote}

In 2004, Vanessa and Caroline traveled to Finland to give Ahlfors'
Fields Medal to the Mathematics Department at the University of
Helsinki.  Six decades after the Medal's audacious exit from Finland,
it had come full circle.  At the suggestion of Olli Lehto, Ahlfors'
longtime colleague and friend, the original Medal was placed in the
Helsinki University Museum, and a replica was made for the Mathematics
Department and is in a glass display case outside the auditorium named
after Lars Ahlfors.

In 2012, Ahlfors' youngest daughter Caroline contacted the Mathematics
Department of LaGuardia Community College about donating her late
father's books.  She was pleased to learn that the books have been
well appreciated by the math faculty, much better than sitting on her
shelf.  She wrote using her husband Frank Mouris' email address.
Later we learned that the couple co-directed \textit{Frank Film},
which won the 1974 Academy Award for Best Animated Short Film.  But
only Frank was recognized; Caroline was overlooked.  Nevertheless,
Ahlfors, with his fame as a Fields Medalist and the author of the
highly-regarded \textit{Complex Analysis}, made this comment: ``Now I
have finally become famous.  I am the father-in-law of an Oscar
winner.''

\section*{Postscript}

After interviewing Vanessa on June 11, 2025, I felt so inspired and I
spent the following days fervently writing the story about Ahlfors'
fascinating family.  On Father's Day (June 15), I shared the
manuscript with Vanessa and expressed my admiration for her beautiful
relationship with Dr.~Ahlfors.  She was touched.

I submitted the article to arXiv, the preprint server managed by
Cornell University, and reached out to mathematicians in Helsinki to
request images of the Fields Medal.  Professor Tuomas Hyt{\"o}nen,
former Head of the Department of Mathematics and Statistics, had been
informed by his predecessor that the Medal was stored in a safe in the
building's basement. However, when they entered the code found in an
envelope, it didn't work.  After replacing the batteries, the safe
remained locked.  Eventually an electrician managed to open the safe,
only to find it empty!

Before calling the police, Professor Hyt{\"o}nen's colleagues
suggested checking with the University Museum.  Riitta-Leena Inki, the
department's science communicator, contacted the Museum and
immediately confirmed that the Medal was in their possession.  Wanting
to verify this firsthand, Professor Hyt{\"o}nen and Ms.~Inki traveled
to the Museum's conservation facilities in Vantaa, about twenty
kilometers from central Helsinki.  There, the curator retrieved the
Medal from storage and handed it to Professor Hyt{\"o}nen, while
Ms.~Inki took several pictures.  For their full account and
photographs of the original Fields Medal, see Inki's article, ``A
Nugget of Gold in Hakkila, Vantaa.''

\renewcommand\refname{Sources}

\subsection*{Photographs}

The photographs shown in this articles were taken by Riitta-Leena
Inki, the science communicator at the University of Helsinki.  She
kindly granted the author the permission to use for this publication.

\newpage
\noindent
\includegraphics[scale=0.4]{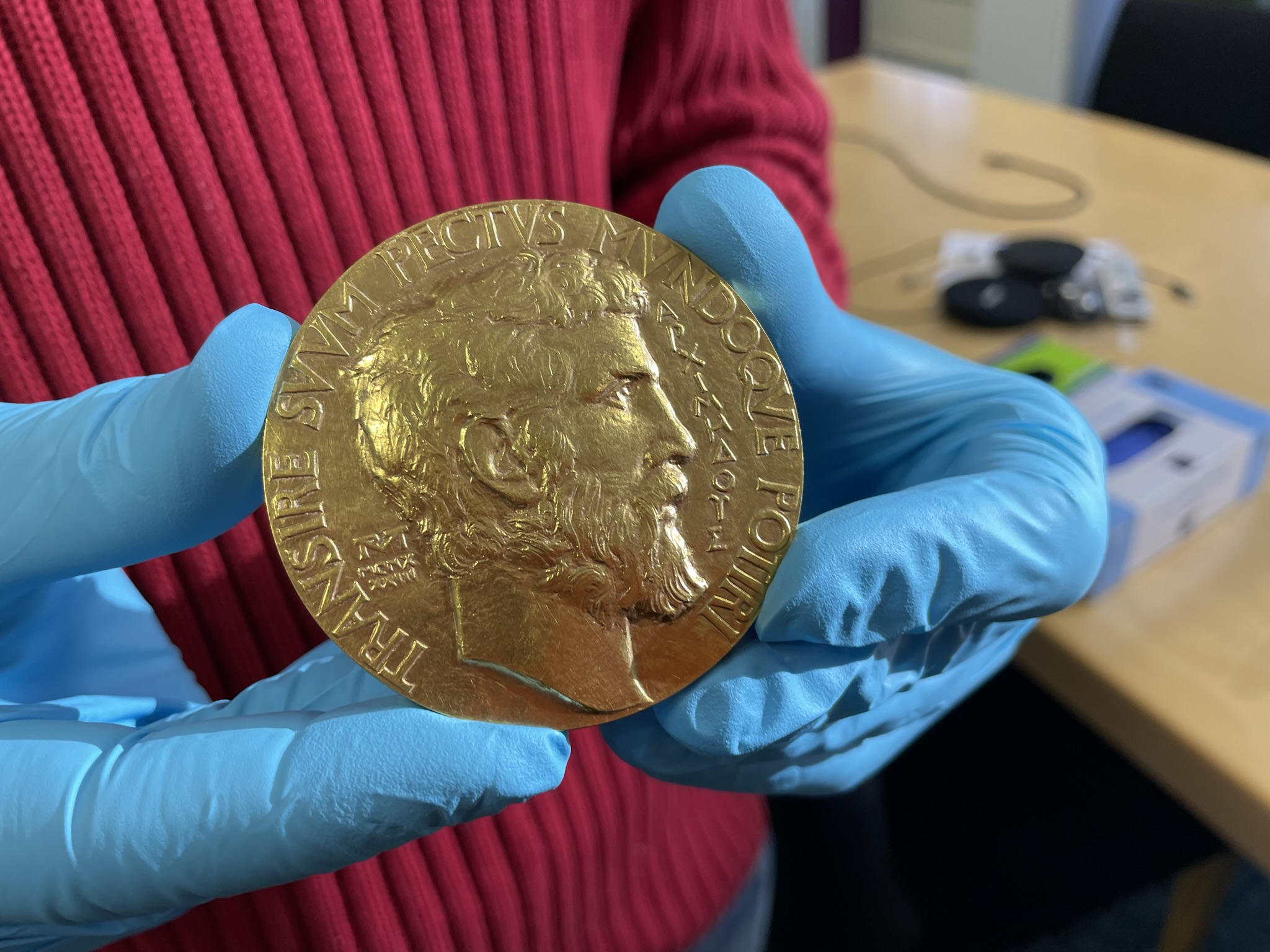}
\includegraphics[scale=0.4]{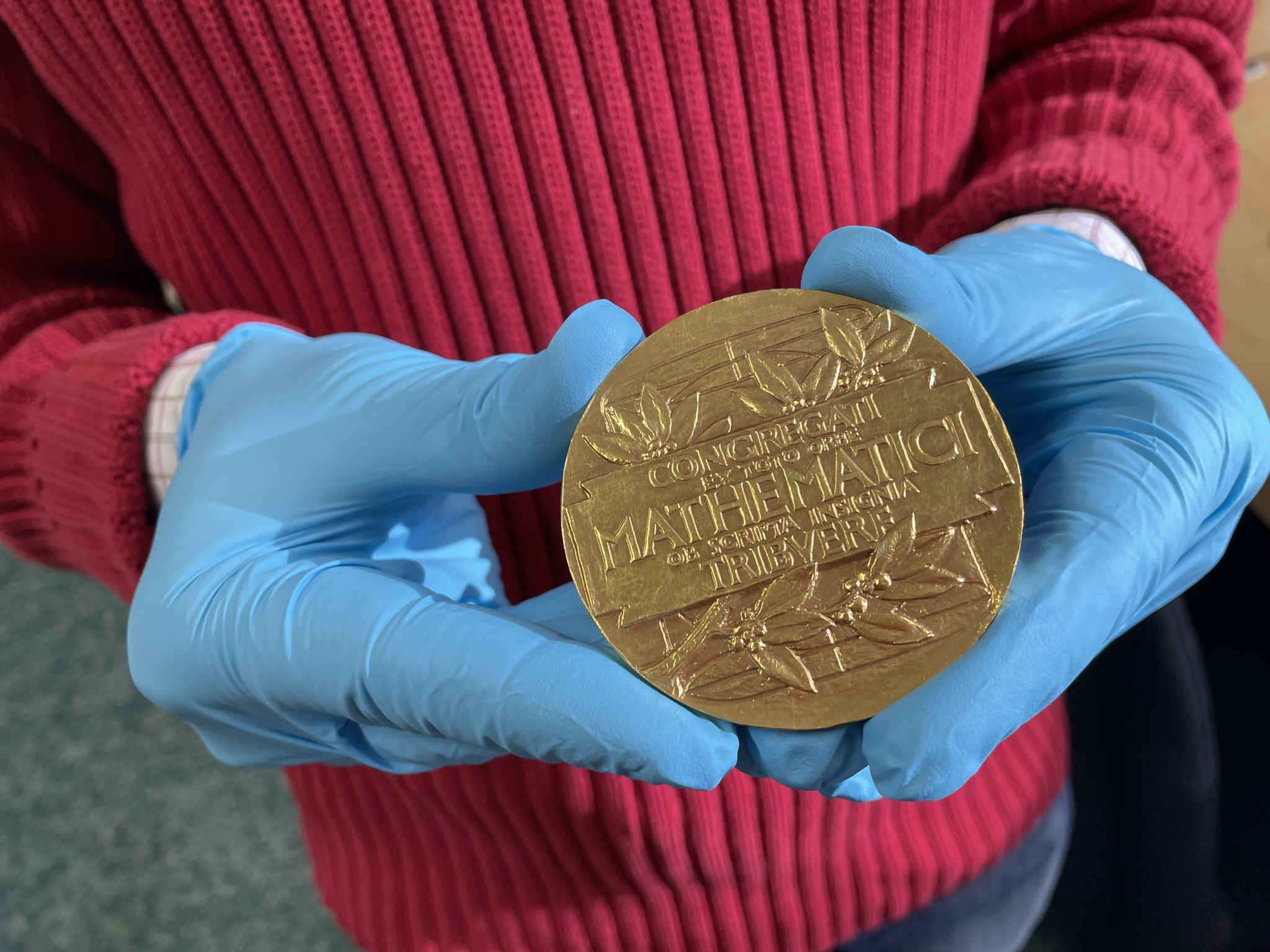}
\includegraphics[scale=0.4]{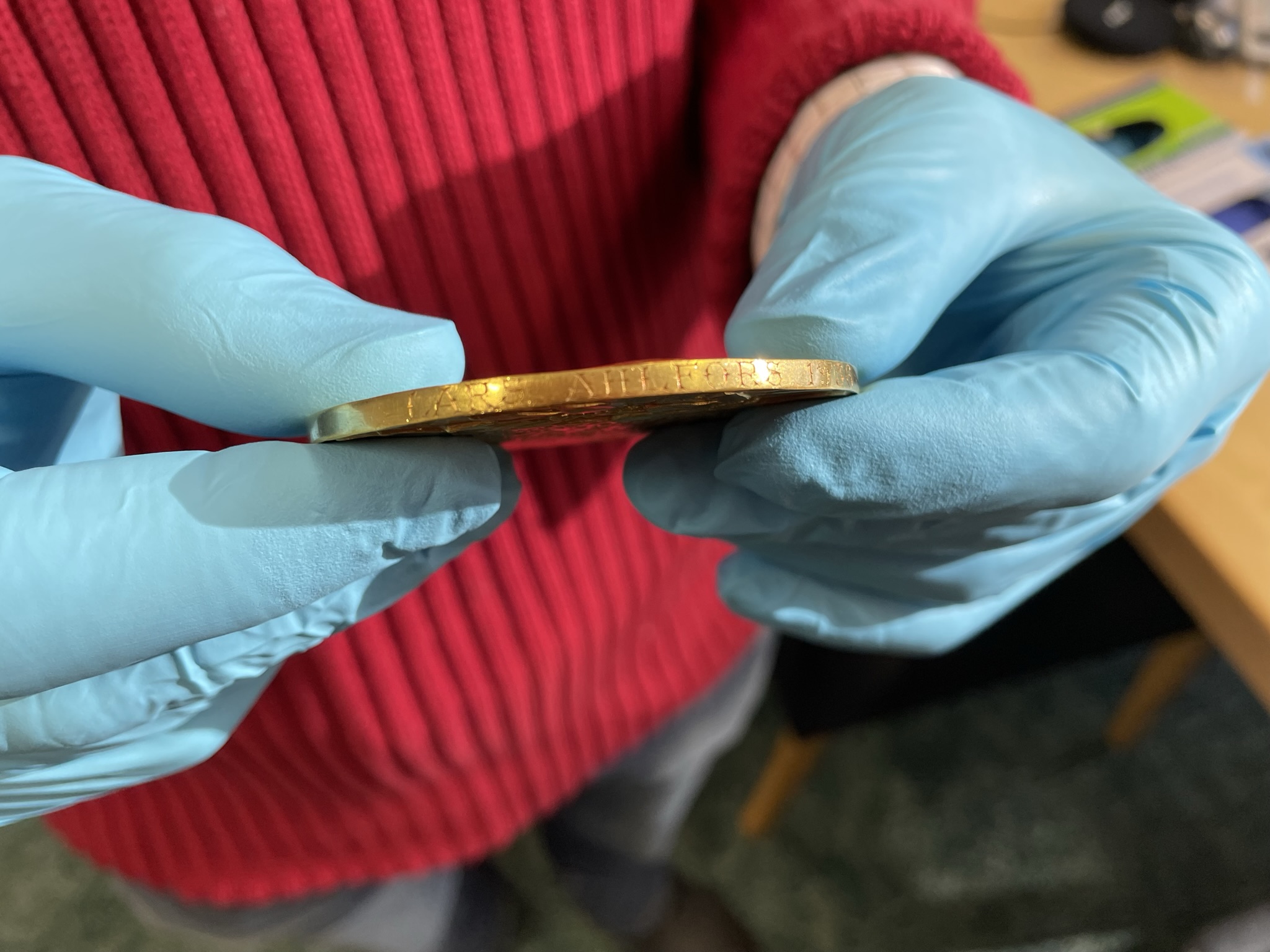}\\
Professor Tuomas Hyt{\"o}nen held the original Fields Medal in his
gloved hands at the conservation facilities of the Helsinki University
Museum in Hakkila, Vantaa.  Photo by Riitta-Leena Inki; courtesy of
Riita-Leena Inki.\\

\noindent
\includegraphics[scale=0.05]{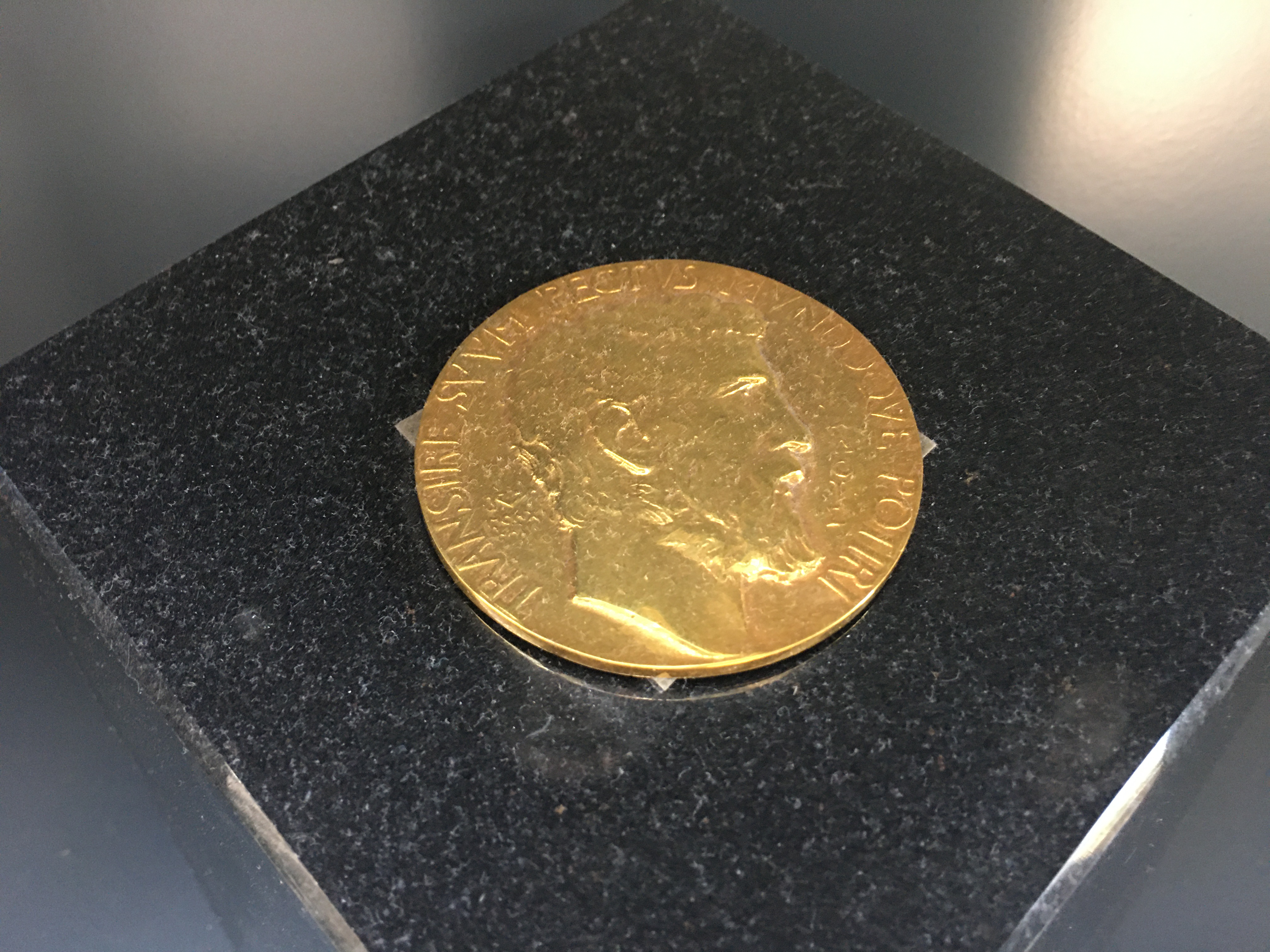}
\includegraphics[scale=0.11]{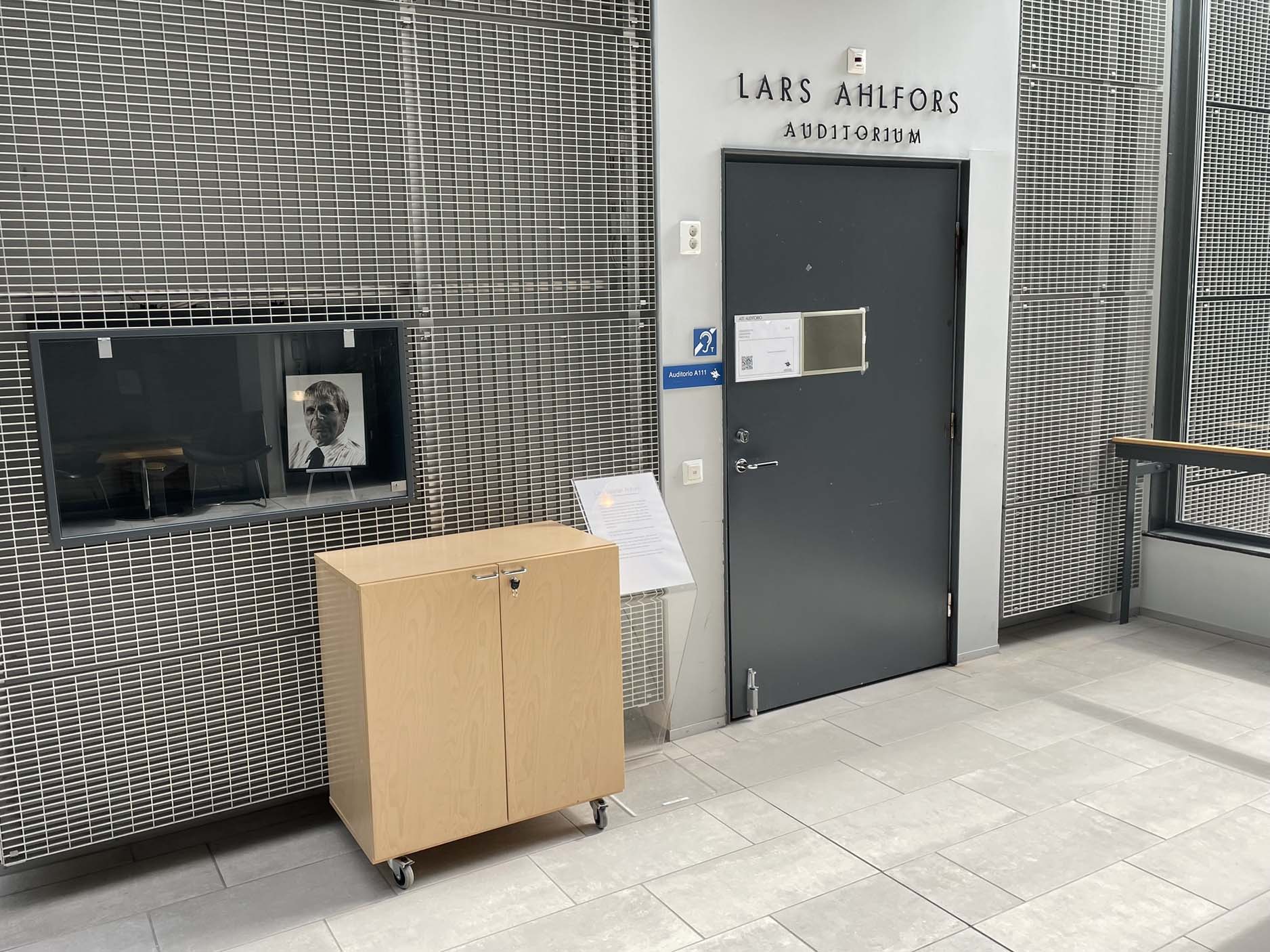}\\
The replica of Lars Ahlfors' Fields Medal is displayed outside the
auditorium named after the mathematician, Exactum A111 at the Kumpula
campus of the University of Helsinki.  Photo by Riitta-Leena Inki;
courtesy of Riita-Leena Inki.

\end{document}